\documentclass[12pt]{amsart}
\usepackage{amsfonts,amssymb,amsmath}
\usepackage{graphics}
\usepackage{epsfig}
\def\pplogo{\vbox{
\halign{##&##\hfil\cr&{
\ppnumber}\cr\rule{0pt}{2.5ex}&\ppdate\cr}
}}
\makeatletter
\def\ps@firstpage{\ps@empty \def\@oddhead{\hss\pplogo}%
  \let\@evenhead\@oddhead 
}
\def\maketitle{\par
 \begingroup
 \def\thefootnote{\fnsymbol{footnote}}
 \def\@makefnmark{\hbox{$^{\@thefnmark}$\hss}}
 \if@twocolumn
 \twocolumn[\@maketitle]
 \else \newpage
 \global\@topnum\z@ 
\quad \vskip.5in
\@maketitle \fi\thispagestyle{firstpage}\@thanks
 \endgroup
 \setcounter{footnote}{0}
 \let\maketitle\relax
 \let\@maketitle\relax
 \gdef\@thanks{}\gdef\@author{}\gdef\@title{}\let\thanks\relax}
\makeatother

\def\ppnumber{\vbox{\baselineskip14pt
\hbox{DUKE-CGTP-00-10}
\hbox{IASSNS-HEP-00/40}
\hbox{math.AG/0007090}
}}
\def\ppdate{July, 2000}

\begin{document}
\title{Geometric Aspects of Mirror Symmetry}
\author{David R. Morrison}

\thanks{Research partially supported by
National Science Foundation grant DMS-9401447 and by the Institute
for Advanced Study, Princeton, New Jersey, USA, whom the author
also thanks for hospitality.}

\address{Center for Geometry and Theoretical Physics, Duke University,
Durham, NC 27708, USA}

\begin{abstract}
The geometric aspects of mirror symmetry are reviewed, with an eye towards 
future developments.  Given a mirror pair $(X,Y)$ of Calabi--Yau threefolds,
the best-understood mirror statements relate certain small corners of the 
moduli spaces of $X$ and of $Y$.  We will indicate how one might go beyond 
such statements, and relate the moduli spaces more globally.  In fact, in 
the boldest version of mirror symmetry (the Strominger--Yau--Zaslow 
conjecture), the Calabi--Yau threefolds $X$ and $Y$ should be directly
related to each other through a very geometric construction. 
\end{abstract}

\maketitle              

\section*{Introduction}

The past twenty years have seen a number of fruitful exchanges of ideas
between pure mathematics and theoretical physics.  On the one hand, deep
results in mathematics arising out of studies in geometry and topology have
had unexpected and powerful applications to theoretical physics.  On the
other hand, insights gleaned from physical models have led to a number of
conceptual revolutions in mathematics, particularly in geometry and
topology.

One of the most exciting of these conceptual revolutions, one which is
still in the process of unfolding, goes by the name ``mirror symmetry.''
Mirror symmetry predicts a completely unexpected relationship between
certain pairs of Calabi--Yau manifolds.\footnote{A {\it Calabi--Yau
manifold}\/ is a K\"ahler manifold with a nowhere-vanishing holomorphic
form of top degree.}  From a mathematical point of view, the connections
between the two ``mirror'' manifolds are extremely indirect, and the effort
to explain them has led to a number of mathematical extensions of the
original mirror phenomenon.  Some of these extensions have proven useful in
physics as well, and physicists have also found more general contexts
in which mirror symmetry can be observed.

Mirror symmetry has become a large field, and there have been a number
of expository surveys 
\cite{MR99k:14063,MR96m:14072,MR99j:32025,MR97j:58013,MR97f:32040,%
Manin,MR97e:32024,MR2000e:14094,MR94c:81194} 
of portions of the theory,
and at least two
books devoted to the subject \cite{MR2000d:14048,MR97i:32026}.
In this paper, we will focus on the geometric connections between mirror pairs,
using the example of the quintic threefold and its mirror as a guide.
We have not attempted to be complete, but have instead concentrated on
a few selected developments which are not covered in other recent
surveys, and which hold much promise for future interesting work.

In the first part of the paper, we review in some detail the known
geometric connections (in the context of toric geometry).  In the second
part of the paper, we outline some new directions in which the geometric
understanding of mirror symmetry has been moving recently, and indicate
some of the important problems remaining to be solved.

Although the issues discussed in this paper are mathematical in nature,
the reader may wish to acquire some knowledge of physics in order to
fully appreciate their context.  A good place to start is 
Deligne et al.\ \cite{MR2000e:81010} (although this assumes some knowledge
of quantum mechanics, which can be obtained, for example, in \cite{sudbery}).
A general introduction to string theory can be found in the
ICM address of Witten \cite{MR89i:81003}, and the standard physics textbooks
on string theory by Green, Schwarz, and Witten
\cite{MR88f:81001} and by
Polchinski \cite{MR99h:81183-4} can be profitably studied.
The more modern aspects of string theory are discussed in detail in \cite{MR98f:81005}.

\section{Duality in Physics and its Mathematical 
   Consequences} \label{sec:duality} 

The origins of mirror symmetry in physics are easy to explain in broad
outline.  Each of the five known superstring theories 
is naturally formulated in ten spacetime
dimensions, and their close cousin ``M-theory'' is naturally formulated as
an eleven-dimensional theory.  When the ten or eleven  dimensions of
the theory take the form of the product of a (small) compact 
Riemannian manifold $X$ with a 
non-compact four-manifold, the geometry of the compact manifold
determines---to first approximation---the physics in four spacetime dimensions.
(This is referred to as ``compactifying the theory'' on $X$.)  The mapping from
geometry to physics need not be one-to-one, however, and this is where the
phenomena of mirror symmetry and other related dualities arise.
(See for example \cite{MR2000c:81225,MR1648148,MR99j:81003,MR99g:81191}.)

Some dualities of this sort involve only a single type of superstring
theory, and are simply statements that the mapping from the geometric
parameter space to the physical parameter space is many-to-one.  For
example, the compactification of the type IIA string theory on a $d$-torus
$X$ has as its natural geometric parameters (in the NS-NS sector) a flat
metric on $X$ and a harmonic element of $H^2(X,U(1))$ known as the
``$B$-field.''  Taking into account the action of $\operatorname{Diff}(X)$,
the geometric parameter space can be described 
in the form $\Gamma_0\backslash\mathcal{D}$, where
\begin{equation}
\mathcal{D}= H^2(X,\mathbb{R})\times GL(d)/O(d)  
\cong O(d,d)/(O(d)\times O(d)) \thinspace ,
\end{equation}
and
\begin{equation}
\Gamma_0=H^2(X,\mathbb{Z})\ltimes GL(d,\mathbb{Z}) \thinspace .
\end{equation}
(This must be modified if $d\ge6$.)  
However, so-called {\em T-dualities}\/ lead to additional 
identifications among points in
this space, and the 
physical NS-NS parameter space
turns out to be
$\Gamma\backslash\mathcal{D}$ instead of $\Gamma_0\backslash\mathcal{D}$, 
where $\Gamma=O(\Lambda^{d,d})$ is the integral
orthogonal group of an even unimodular lattice of signature $(d,d)$.
(Intrinsically, 
$\Lambda^{d,d}\cong H^1(X,\mathbb{Z})\oplus H_1(X,\mathbb{Z})$.)

Other dualities involve two types of superstring theories.  In fact,
another version of 
T-duality is also of this sort: in its most primitive form,
T-duality gives an identification between the type IIA string theory
compactified on a $d$-torus and the type IIB string theory compactified on
a (geometrically different) $d$-torus, when $d\ge1$.\footnote{There are 
also versions of T-duality
which similarly relate the two kinds of heterotic string theory.}
In the case of a rectangular metric on the torus (and a $B$-field value
of zero), the statement is  that a rectangular torus with
radii $(r_1,\dots,r_d)$ is dual to one with radii $(1/r_1,\dots,1/r_k,r_{k+1},
\dots,r_d)$, where the duality interchanges IIA and IIB if $k$ is odd, and
preserves both IIA and IIB theories if $k$ is even.

The geometric spaces which appear in the duality known as mirror symmetry
 are Riemannian manifolds of dimension $2n$ with holonomy contained in
 $SU(n)$, often called {\it Calabi--Yau manifolds}.  (The theorems of
Calabi \cite{calabi} and Yau \cite{yau} lead to an alternate
characterization of these as K\"ahler manifolds with a nowhere vanishing
holomorphic $n$-form.)  When $n$ is odd, mirror symmetry relates the type
IIA string theory compactified on one Calabi--Yau manifold $X$ to the type
IIB string theory compactified on another Calabi--Yau manifold $Y$, and
vice versa.  When $n$ is even, the IIA theories on $X$ and on $Y$ are
related, as are the IIB theories on $X$ and on $Y$.  Note that the
``other'' manifold $Y$ will in some cases be topologically the same as
$X$, but typically the geometric parameters (a metric and some harmonic
$p$-forms) will be different.

In the case of K3 surfaces, mirror symmetry becomes an assertion that the
physical parameter space is obtained from the geometric one by additional
identifications \cite{MR97i:81128}.  For the type IIA theory, the geometric
parameter space is again described in terms of a metric on the K3 surface
$X$ and a $B$-field.  Each metric is characterized by the corresponding
$3$-dimensional plane of self-dual harmonic $2$-forms inside the
$22$-dimensional space $H^2(X,\mathbb{R})$.  Including the limiting
``orbifold metrics,'' the geometric parameter space is then described as
$\Gamma_0\backslash\mathcal{D}$, where this time
\begin{equation}
\begin{aligned}
\mathcal{D}&=H^2(X,\mathbb{R})\times\mathbb{R}^+\times  
O(3,19)/(O(3)\times O(19)) \\ 
&\cong O(4,20)/(O(4)\times O(20))
\end{aligned}
\end{equation}
and
\begin{equation}
\Gamma_0=H^2(X,\mathbb{Z})\ltimes O^+(H^2(X,\mathbb{Z})) \thinspace .
\end{equation}
(Here $O^+(H^2(X,\mathbb{Z}))$ is an index $2$ subgroup of
$O(H^2(X,\mathbb{Z}))$, isomorphic to
$\operatorname{Diff}(X)/\operatorname{Diff}^0(X)$
\cite{Mat:d,Borcea:d,Don:inv}.)  The total cohomology $H^*(X,\mathbb{Z})$
has the structure of an even unimodular lattice of signature $(4,20)$
\cite{Muk:bundles}, and the physical parameter space becomes
$\Gamma\backslash\mathcal{D}$, with $\Gamma=O(H^*(X,\mathbb{Z}))$.  This
physical mirror symmetry for K3 surfaces has many beautiful mathematical
consquences, including assertions that families of algebraic K3 surfaces
defined by conditions on the Picard group have ``mirror dual families'',
and a connection to Arnold's ``strange duality''
\cite{MR97i:81128,MR97i:14024,Kobayashi}.

In general,
the assertion that two Calabi--Yau manifolds $X$ and $Y$ (with their geometric
parameters appropriately specified) form a ``mirror pair'' is 
really a statement
about physics rather than about mathematics.  The geometric data we have used
above is a first approximation to the true physical parameters, and at present
it is not known how to completely specify the physical parameter space, or 
even to fully
test a proposed equivalence between two physical theories.

However, there are a number of purely mathematical statements which can be
extracted from what we know of the physics of a mirror pair, and these 
statements
can be tested in examples, and studied in their own right.  Several such 
mathematical spinoffs of mirror symmetry have
been the subject of intense research during the past decade.

The reason mirror symmetry has attracted so much attention from mathematicians
is its truly surprising nature.  Prior to the discovery of this phenomenon
by physicists \cite{Dixon,LVW,MR91m:14062,MR91h:32018}, it was not even 
suspected that most Calabi--Yau manifolds would
come in closely-related pairs.  Moreover, some of the consquences of the
physical mirror equivalence were shocking: a generating function for the
number of holomorphic curves of various fixed degrees on one of the Calabi--Yau
manifolds could be computed using the integrals of a top-degree holomorphic
form on the mirror partner---a totally unexpected result.

\section{The Quintic Threefold and its Mirror} \label{sec:quintic}

Among the simplest known examples of Calabi--Yau threefolds are the
{\it quintic threefolds}: hypersurfaces of degree $5$ in complex
projective space $\mathbb{P}^4$.  The complex projective space, including
a choice of K\"ahler metric on it, can be described by means of
{\it symplectic reduction}: begin with the effective group action of $U(1)$ on
$\mathbb{C}^5$ defined by
\begin{equation}\label{eq:groupaction}
e^{i\theta}: (x_1,\dots,x_5)\mapsto(e^{i\theta}x_1,
\dots,e^{i\theta}x_5)
\end{equation}
which admits a moment map 
$\mu:\mathbb{C}^5\to\mathfrak{g}^*\cong\mathbb{R}^1$
given by
\begin{equation}
\mu(x_1,\dots,x_5)=\frac12\sum_{j=1}^5|x_j|^2 \thinspace ;
\end{equation}
then $\mathbb{P}^4=\mu^{-1}(r)/U(1)$ with K\"ahler form
induced by the canonical K\"ahler form $\sum dz_i\wedge d\overline{z}_i$
on $\mathbb{C}^5$.  The image of the moment map contains all positive
values of $r$, and all K\"ahler classes on $\mathbb{P}^4$ are produced
in this way.
The quotient space is a {\it toric variety}: it admits an action
of a complex torus $T\cong (\mathbb{C})^*{}^4$, the complexification
of $U(1)^5/U(1)$, with an open dense orbit.

The group action \eqref{eq:groupaction} which was used to make this 
construction can be
specified by means of an injective homomorphism
$\alpha:\mathbb{Z}\to\mathbb{Z}^5$ (which is instrinsically described
as $\pi_1(U(1))\to\pi_1(U(1)^5)$),
or equivalently, by the cokernel $\mathbb{Z}^4=\pi_1(T)$ of
$\alpha$, together with the images in $\pi_1(T)$
of the standard basis vectors of $\mathbb{Z}^5$.
Concretely, $\alpha(1)=(1,1,1,1,1)$, and
we may choose a basis of $\pi_1(T)$ so
that the image vectors are given by
\begin{equation}\label{eq:vec}
\begin{gathered}
(1,0,0,0),\quad (0,1,0,0),\quad (0,0,1,0), 
\\
(0,0,0,1),\quad (-1,-1,-1,-1) \thinspace .
\end{gathered}
\end{equation}
In fact (as we shall see in Section \ref{sec:next}), in toric geometry it 
is the convex polytope spanned by those image vectors which plays a 
crucial r\^ole.

The Calabi--Yau hypersurface is specified by a homogeneous polynomial
$f(x_1,x_2,x_3,x_4,x_5)$ of degree $5$. The possible monomials in that
polynomial form another convex polytope, the {\it Newton polytope}\/
of the hypersurface, whose vertices correspond to the monomials
$x_1^5$, \dots, $x_5^5$.  It is convenient to select one monomial
such as $x_1x_2x_3x_4x_5$ and study the quotient
\begin{equation}\label{eq:quot}
\frac{f(x_1,x_2,x_3,x_4,x_5)}{x_1x_2x_3x_4x_5} \thinspace ,
\end{equation}
since this quotient is invarient under $U(1)$.
In fact, each Laurent monomial occurring in \eqref{eq:quot}
is a character of $T$ and can be labelled by the
corresponding element of 
$\operatorname{Hom}(T,\mathbb{C}^*)\cong
\operatorname{Hom}(\pi_1(T),\mathbb{Z})$.
In the basis dual to the basis of $\pi_1(T)$ used in \eqref{eq:vec},
the Newton polytope is spanned by
\begin{equation}\label{eq:vecdual}
\begin{gathered}
(4,-1,-1,-1),\quad (-1,4,-1,-1),\quad (-1,-1,4,-1),
\\
(-1,-1,-1,4),\quad (-1,-1,-1,-1) \thinspace ,
\end{gathered}
\end{equation}
but of course also contains many other lattice points corresponding to
other monomials.

To form the mirror partner \cite{MR95c:14046}, we reverse the r\^oles 
of the polytopes
spanned by \eqref{eq:vec} and \eqref{eq:vecdual}.  That is, we
use the data in \eqref{eq:vecdual} to specify the
symplectic reduction, and we use the data in \eqref{eq:vec} to
determine a hypersurface.  There are some subtleties
associated with selecting $r$, the point in the image of the
moment map, to which we return in Section \ref{sec:next}.
In brief, this toric variety turns out to be\footnote{Actually, the
mirror of the quintic was originally constructed as a quotient,
rather than as a toric hypersurface \cite{MR91h:32018}.} a resolution of
singularities of $\mathbb{P}^4/(\mathbb{Z}_5)^3$, and
the mirror hypersurface is the image of a hypersurface with
equation
\begin{equation} \label{eq:mirhyp}
c_0x_1x_2x_3x_4x_5+c_1x_1^5+\cdots+c_5x_5^5=0 \thinspace .
\end{equation}
(We will explain how this is derived in Section~\ref{sec:explain} below.)

If we let $(\mathbb{C}^*)^5$ act on $x_1, \dots, x_5$, the set of
hypersurfaces of the form \eqref{eq:mirhyp} is permuted.  So we can regard
$(\mathbb{C}^*)^5$ as acting on the set of coefficients
$\{(c_0,c_1,\dots,c_5)\}$; the quotient by this action is the moduli space
of complex structures.\footnote{There are a number of technical
difficulties involved in forming this quotient, which we do not discuss
here; see \cite{MR95m:14012,MR95b:14029,MR2000d:14048}.}

\begin{figure}
  \centerline{\epsfxsize=4cm\epsfbox{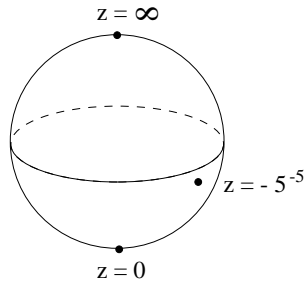}}
  \caption{The moduli space of complex structures on the quintic mirror.}
  \label{fig:moduli}
\end{figure}

In our example, the compactified moduli space is isomorphic to
$\mathbb{P}^1$, and a coordinate on it is provided by
\begin{equation}
z=\frac{c_1c_2c_3c_4c_5}{c_0^5} \thinspace ,
\end{equation}
which generates the $(\mathbb{C}^*)^5$-invariant functions on
$\{(c_0,c_1,\dots,c_5)\}$.  The mirror hypersurface acquires additional
singularities at $z=0$ and at $z=-5^{-5}$; moreover, the $z=\infty$ limit
coincides with the Fermat hypersurface.  The moduli space is illustrated in
Figure~\ref{fig:moduli}.

Now that we have described a mirror pair of Calabi--Yau manifolds, we
can ask:  what kind of properties are
they supposed to share (according to predictions from physics), and which
ones can we verify?

The first properties have to do with the Hodge numbers of the Calabi--Yau
threefolds.  If $X$ and $Y$ are mirror partners, then we expect that
$h^{1,1}(X)=h^{2,1}(Y)$ and $h^{2,1}(X)=h^{1,1}(Y)$.  
This property holds for our example mirror pair.

In fact, $h^{1,1}(X)$ is the dimension of the K\"ahler moduli space of
$X$, while $h^{2,1}(X)$ is the dimension of its complex structure moduli space.
The 
equality of Hodge numbers
is just one part of a stronger property, which
asserts that these moduli spaces for $X$ and for $Y$ should coincide in 
some appropriate sense,
with the r\^oles of K\"ahler and complex structure moduli reversed.

Thanks to the expectation in physics that the geometric version of the
(complexified) K\"ahler moduli space requires so-called quantum corrections
which may affect its structure, the best prediction which can be made is
that the moduli spaces should agree locally in the neighborhood of some
boundary point.  That is, for an appropriate compactification of
the moduli space of complex strutures of $Y$, there should be a boundary point
with a neighborhood which coincides
with a neighborhood in the complexified K\"ahler moduli space of $X$.
In the case of the mirror quintic, this ``large complex structure limit
point'' is at $z=0$.  (The singularity at $z=-5^{-5}$ is of a different
nature.) 

\section{Periods and Hypergeometric Systems}

One of the key features of the complex structure moduli space is the behavior
of the period integrals on the Calabi--Yau manifold as a function of 
parameters.
These period integrals always satisfy an algebraic differential equation,
and in the case of toric complete intersections, these algebraic differential
equations turn out to be hypergeometric systems \cite{MR94m:14067}, in
fact, the 
generalized hypergeometric systems which have been studied in detail by
Gel'fand, Zelevinski\u\i, and Kapranov \cite{GZK}.

For example, in the case of the mirror quintic, the differential equation
satisfied by the periods $\Phi(z)$ can be written as $\mathcal{D}\Phi=0$ where
\begin{equation}
\mathcal{D}=\left(z\frac d{dz}\right)^4
-5z\left(5z\frac d{dz}+1\right)
\left(5z\frac d{dz}+2\right)
\left(5z\frac d{dz}+3\right)
\left(5z\frac d{dz}+4\right) \thinspace .
\end{equation}

It is easy to find a single power series solution near $z=0$
\cite{MR96g:32037}: 
\begin{equation}
\Phi(z)=\sum_{n=0}^\infty \frac{(5n)!}{(n!)^5}\,z^n
\end{equation}
but the other three solutions are elusive.  The recursion relations
implied by the equation lead one to a formal power series of the form
\begin{equation}\label{eq:formal}
\Phi_\alpha(z)=\sum_{n=0}^\infty \frac{(5\alpha +1)(5\alpha +2)
\cdots(5\alpha +5n)}{[(\alpha +1)(\alpha +2)\cdots(\alpha +n)]^5}
\,z^{\alpha +n} \thinspace ;
\end{equation}
one finds that $\mathcal{D}(\Phi_\alpha(z))=\alpha ^4z^\alpha $ and so we 
must have $\alpha ^4=0$
in order to obtain a solution.  In fact \cite{MR2000b:14013}, the formal 
solution
\eqref{eq:formal} can be interpreted with $\alpha $ taken from the ring
$\mathbb{C}[\alpha ]/(\alpha ^4)$ as follows:  each coefficient 
\begin{equation}
\frac{(5\alpha +1)(5\alpha +2)\cdots(5\alpha +5n)}{[(\alpha +1)
(\alpha +2)\cdots(\alpha +n)]^5}
\end{equation}
can be evaluated in that
ring, and written as a polynomial in $\alpha $ of degree $3$.  This gives
a four-dimensional space of multi-valued solutions.  The multi-valuedness
arises due to the interpretation of $z^\alpha $ as $\exp(\alpha \log z)$.
 
In other words, if we replace $\log z$ by $\log z + 2\pi i$, then the
solution $\Phi_\alpha(z)$ is multiplied by $e^{2\pi i\alpha}$.  
That is, the monodromy on the set of
solutions near $z=0$ is represented by the action (by multiplication) of
$e^{2\pi i\alpha}$ on the ring $\mathbb{C}[\alpha]/(\alpha^4)$.

Mirror symmetry proposes an identification between the cohomology ring of
one Calabi--Yau manifold (in our example, the quintic), and the ring
structure on the period integrals of the other.  There are factors of $2\pi i$
which must be introduced when comparing these rings: the generators
$\alpha_j$ which
lead to monodromy transformations of the form $e^{2\pi i\alpha_j}$
must be mapped to $\lambda_j/2\pi i$, where the $\lambda_j$ are generators of
the integral cohomology of the other manifold.
Since the cohomology ring
of the quintic takes the form $\mathbb{C}[\lambda]/(\lambda^4)$, 
the proposed identification
 is verified in our example, with $\lambda=2\pi i\alpha$.
  It is known to hold for a broad class
of Calabi--Yau complete intersections
\cite{MR91g:81107,MR95i:32028,MR96d:32028,MR2000b:14013}. 

There is no nice form for a power series solution near $z=-5^{-5}$.
However, we {\it can}\/ find explict power series solutions to 
$\mathcal{D}\Phi=0$
near $z=\infty$.  Begin with a power
series of the form
\begin{equation}
\Phi_\alpha(z)=\sum_{m=0}^\infty \frac{[(\alpha )(\alpha +1)
\cdots(\alpha +m-1)]^5}{(5\alpha )(5\alpha +1)\cdots(5\alpha +5m-1)}
\,z^{-\alpha -m} \thinspace ,
\end{equation}
and calculate 
\begin{equation}
\mathcal{D}(\Phi_\alpha(z))=-5(-5\alpha +1)(-5\alpha +2)
(-5\alpha +3)(-5\alpha +4)z^{-\alpha +1} \thinspace .
\end{equation}
Thus, to get a solution, we need
\begin{equation}
(5\alpha -1)(5\alpha -2)(5\alpha -3)(5\alpha -4)=0 \thinspace ,
\end{equation}
and we get four independent solutions this way.  These are again multi-valued,
since a fifth root of $z$ must be extracted to define the solution.
The monodromy this time is given by the action of $e^{2\pi i\alpha}$ on the
ring
\begin{equation}
\mathbb{C}[a]/(5\alpha -1)(5\alpha -2)(5\alpha -3)(5\alpha -4) \thinspace ;
\end{equation}
this gives a monodromy transformation of order $5$.

At all other points in the moduli space there are four independent solutions
to $\mathcal{D}\Phi=0$.  Thus, it is only at the three points marked in 
Figure~\ref{fig:moduli} that there is a monodromy transformation associated
to a loop around the point.

Among the quantities of importance in physics are the so-called 
topological correlation functions,
which can be calculated from the periods near $z=0$.  
A close study of these functions led to the striking
predictions about enumerative geometry which first awakened the interest of
many mathematicians in this field.  (See \cite{MR93b:32029} for 
the original calculation, and 
  \cite{MR93j:14047,MR96g:32037,MR98a:14055} for accounts of these predictions
intended for mathematicians.)
These predictions were subsequently
verified in the celebrated Mirror Theorem of  Givental \cite{MR97e:14015} and 
 Lian--Liu--Yau \cite{MR99e:14062}.
The proofs are actually logically independent of mirror
symmetry---they consist in showing that the enumerative geometry of rational
curves on the quintic (or on other Calabi--Yau manifolds) has a rich
structure which reproduces the predictions made by mirror symmetry.
In this way, one verifies that there are 2875 lines on the general
quintic threefold, 609250 conics, 317206375 twisted cubics, and so on.

For the mirror partner, the enumerative geometry (or the topological correlation
functions) can be encoded in a structure known as the {\it quantum
cohomology ring}\/ of the Calabi--Yau 
manifold.\footnote{For an introduction to quantum cohomology, 
see
\cite{MR98h:58031}.}  (In fact, the quantum
cohomology is not only a ring, but also a graded Frobenius algebra.)
The quantum cohomology is a deformation of the usual cohomology ring,
with the parameter space being a (formal) neighborhood of the large
complex structure limit point in the moduli space.  The counts of numbers
of rational curves govern the coefficients in the deformation,
as well as the ``trace'' in the Frobenius algebra structure.

\section{Triangulations and Flops} \label{sec:next}

Mirror symmetry is expected to be a two-way street, so we should be able
to reverse the r\^oles of the quintic hypersurface and its mirror partner.
However, things immediately become much more complicated.

\begin{figure}
\setlength{\unitlength}{.38em}

\begin{center}
\begin{picture}(40,30)

\thinlines

\put(0,0){\line(1,0){40}}
\put(40,0){\line(-2,3){20}}
\put(20,30){\line(-2,-3){20}}

\put(4,6){\line(1,0){32}}
\put(32,0){\line(-2,3){16}}
\put(24,24){\line(-2,-3){16}}

\put(0,0){\line(2,1){12}}
\put(40,0){\line(-2,1){12}}
\put(20,30){\line(0,-1){12}}

\put(16,0){\line(2,3){12}}

\put(20,6){\line(-2,-1){12}}
\put(20,6){\line(2,-1){12}}

\put(32,12){\line(-1,0){24}}

\put(24,12){\line(2,-1){12}}
\put(24,12){\line(0,1){12}}

\put(12,18){\line(2,-3){12}}

\put(16,12){\line(0,1){12}}
\put(16,12){\line(-2,-1){12}}

\end{picture}
\hfill
\begin{picture}(40,30)

\thinlines

\put(0,0){\line(1,0){40}}
\put(40,0){\line(-2,3){20}}
\put(20,30){\line(-2,-3){20}}

\put(4,6){\line(1,0){32}}
\put(32,0){\line(-2,3){16}}
\put(24,24){\line(-2,-3){16}}

\put(8,0){\line(-2,3){4}}
\put(32,0){\line(2,3){4}}
\put(16,24){\line(1,0){8}}

\put(16,0){\line(2,3){12}}

\put(16,0){\line(-2,3){4}}
\put(24,0){\line(2,3){4}}

\put(32,12){\line(-1,0){24}}

\put(32,12){\line(-2,-3){4}}
\put(20,18){\line(1,0){8}}

\put(12,18){\line(2,-3){12}}

\put(12,18){\line(1,0){8}}
\put(8,12){\line(2,-3){4}}

\end{picture}
\end{center}
\caption{Parts of two possible triangulations for the mirror of the 
   quintic}\label{fig1}
\end{figure}
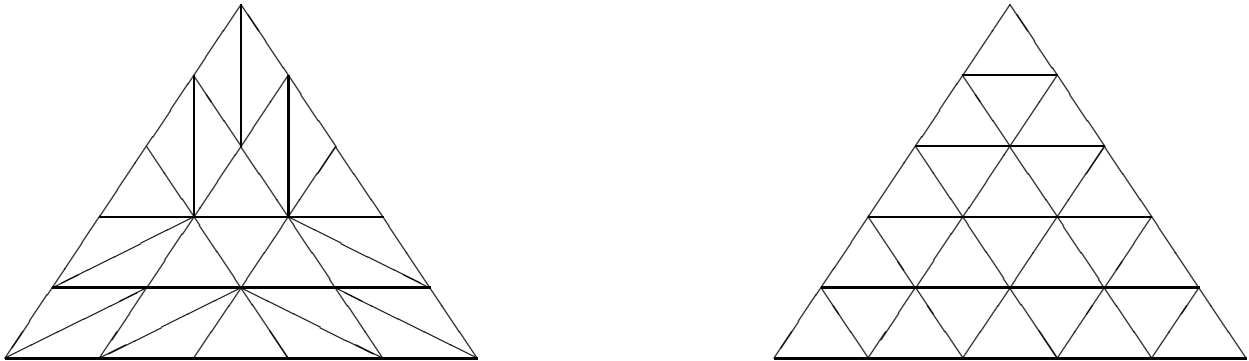

As mentioned in Section \ref{sec:quintic}, the 
mirror partner is a resolution of singularities of a hypersurface in
$\mathbb{P}^4/(\mathbb{Z}_5)^3$.  There are many many possible choices
of resolution, corresponding (according to toric geometry)
to the possible triangulations of the
polytope spanned by \eqref{eq:vecdual}.
One possible choice, described explicitly in \cite{MR93j:14047},
is depicted on the left side of Figure~\ref{fig1}; another, perhaps
more symmetric choice, is depicted on the right side.  Computing
all possible triangulations would appear to be beyond the range of current computers
\cite{MR2000d:14048}.

For each possible choice of triangulation, there is an asoociated set of
moment map values $r$ which give rise to that triangulation.  We thus get
a vast number of possible Calabi--Yau manifolds and associated K\"ahler
moduli spaces, each of them a mirror of the original quintic threefold.
However, all of these Calabi--Yau manifolds are related by flops (the
simplest type of birational operation), and it is understood in physics
that the corresponding K\"ahler moduli spaces can be sewn together into
a larger space which contains all of them \cite{MR95a:81261,MR95i:32027} 
(see also  \cite{MR97a:32020} for a more mathematical account).

In fact, this is very natural from the point of view of mirror symmetry.
In the case of the quintic, for example, there is a complex structure
moduli space of dimension $101$ with a vast number of ``large complex
structure limit'' boundary points.  To each of these is associated a
possible birational model of the mirror.

The enumerative geometry implications of mirror symmetry could be
calculated for each of these mirrors---each would involve a calculation
which involves only a small neighborhood of the corresponding large complex
structure limit point.  Fascinating though such computations, are,
our focus here will be on ways to relate the {\it entire}\/ moduli
spaces.

Some preliminary hints that the entire moduli spaces are related to each other 
 can be found by studying the combinatorics of the moduli space.
In the case of the quintic, or more generally for complete intersection
in toric varieites, the moduli space is itself a toric variety,
associated to a combinatorial structure known as the {\it secondary polytope}\/
or {\it secondary fan}\/ \cite{MR95e:14045}.  This structure can be
studied from the point of view of a polynomial whose coefficients
help to describe the fan (really, it is the discriminant locus
of this polynomial which is relevant), or dually, from the point of view of
enumerating possible triangulations of a fixed polytope (such as 
\eqref{eq:vecdual}).  The two points of view correspond to the two
sides of the mirror picture we have presented.

\section{The Gauged Linear Sigma Model} \label{sec:explain}

The description we have given of the quintic and its mirror is a special
case of a construction of Batyrev \cite{MR95c:14046} giving mirror partners
for hypersurfaces in toric varieties.  There is a natural generalization
of this, due to Batyrev and Borisov \cite{MR98b:14033}, to the case of
complete intersections in toric varieties.  The generalization is also
closely related to an important construction in physics known as the
gauged linear sigma model \cite{MR95a:81261}.

We illustrate the basic construction with the case of the quintic
hypersurface.  The homogeneous polynomial $f(x_1,\dots,x_5)$ can be
regarded either as describing a section of the line bundle
$\mathcal{O}(-K_{\mathbb{P}^4})$, or as defining a complex-valued function
\begin{equation}
W(x_0,x_1,\dots,x_5):=x_0f(x_1,\dots,x_5)
\end{equation}
on the total space of the line bundle $\mathcal{O}(K_{\mathbb{P}^4})$ 
(with fiber
coordinate $x_0$).
In the latter interpretation, the
Calabi--Yau threefold coincides with the critical set $\operatorname{Crit}(W)$
of the function $W$.  In fact, if the polynomial $f$ is
transverse,\footnote{That is, if the only common zero of its partial
derivatives is the origin.} the critical set away from the origin in
$\mathbb{C}^5$ is defined by $x_0=f(x_1,\dots,x_5)=0$. 

(A similar construction leads from a complete intersection
$f_0=\cdots=f_{k-1}=0$ in a projective toric variety, with $f_j$ a section of
$\mathcal{O}(L_j)$, to the polynomial function
$W=x_0f_0+\cdots+x_{x-1}f_{k-1}$ on the total space of the bundle
$\mathcal{O}(-L_0)\oplus\cdots\oplus\mathcal{O}(-L_{k-1})$.  The critical
set $\operatorname{Crit}(W)$ again coincides with the original complete intersection.)

The ambient space $\mathcal{O}(K_{\mathbb{P}^4})$ can itself be described by a
quotient construction.  We again describe it  using symplectic reduction:
there is an action of $U(1)$ on $\mathbb{C}^6$ defined by
\begin{equation}
e^{i\theta}: (x_0,x_1,\dots,x_5)\mapsto(e^{-5i\theta}x_0,e^{i\theta}x_1,
\dots,e^{i\theta}x_5)
\end{equation}
which admits a moment map 
$\mu:\mathbb{C}^6\to\mathfrak{g}^*\cong\mathbb{R}^1$
given by
\begin{equation}
\mu(x_0,x_1,\dots,x_5)=\frac12\big(-5|x_0|^2+\sum_{j=1}^5|x_j|^2\big)
\thinspace , 
\end{equation}
and $\mathcal{O}(K_{\mathbb{P}^4})=\mu^{-1}(r)/G$ for appropriate values 
of $r$.
Note that the polynomial function $W:\mathbb{C}^6\to\mathbb{C}$ is 
$G$-invariant.

The general version of this construction goes as follows: given a subgroup $G$
of $U(1)^n$, a point $r\in\mathfrak{g}^*$ in the image of the moment map,
 and a $G$-invariant polynomial $W$ defining a function on $\mathbb{C}^n$,
we can study $\operatorname{Crit}(W)$ on the quotient variety
$\mu^{-1}(r)/G$.  
If the polynomial $W$ is given explicitly as
\begin{equation} \label{eq:W}
W(x_0,\dots,x_{n-1})=\sum_{j=0}^{m-1}c_j\prod_{k=0}^{n-1}x_k^{p_{jk}}
\end{equation}
(with $c_j\ne0$), then
the combinatorics in this construction are 
essentially
captured by the $m\times n$ matrix of exponents $P=(p_{jk})$. 
(We assume that $W$ contains
enough monomials so that the rank of $P$ is $d:=n-\dim G$.)
As in Section \ref{sec:quintic}, the coefficients $c_j$ are somewhat
redundant:  there is a group which acts on the space of polynomials of 
the form \eqref{eq:W}, and we must form the quotient by this group action.
The true coordinates on the complex structure moduli space are the
invariant quantities $z_\ell$ for this group action.

The conditions which this data must
satisfy are that the monomials in $W$ generate a Gorenstein cone, and that
the dual of this cone also be Gorenstein.
In terms of the matrix $P$, this means that there must exist rational vectors
$\mu$ and $\nu$ such that 
$P\mu={}^t(1,\dots,1)$, and ${}^t\!\nu P=(1,\dots,1)$.
It turns out that these conditions imply that whenever
 $\operatorname{Crit}(W)$ is a manifold of dimension
$d-2 ({}^t\!\nu P\mu)$, it is a Calabi--Yau manifold.
The triangulations of a fixed polytope which led to the ``secondary fan''
description in the previous section are now replaced by triangulations
of the Gorenstein cone.

To encode the group $G$ in the same combinatorial structure, we introduce a
basis $x^{t_\alpha}$ of $G$-invariant Laurent monomials on $\mathbb{C}^n$;
then we can write
\begin{equation}
x^{p_j}=\prod_{\alpha=1}^d (x^{t_\alpha})^{s_{j\alpha}}
\end{equation}
for each of the monomials occurring in $W$.  This gives a factorization
of $P$ as the product of an $m\times d$ matrix $S$ and a $d\times n$
matrix $T$, with the
group $G$ being completely determined by $T$.  Changing the basis of Laurent
monomials alters $(S,T)$ to $(SL^{-1},LT)$ for some $L\in GL(d,\mathbb{Z})$.

The gauged linear sigma model \cite{MR95a:81261} is a physical theory built
from the group $G$ and its action on the $x$'s.  It is expected that at
low energies, this theory will agree with (the perturbative part of)
type II string theory compactified on the associated Calabi--Yau manifold.

There is a mirror partner of a gauged linear sigma model, whose construction
is essentially due to Batyrev and Borisov \cite{MR98b:14033} (see also 
\cite{MR96a:32043,MR97a:32030,MR98i:81236}).  To describe 
the mirror partner, one merely replaces $P$, $S$, and $T$ by their 
transposes.  
The dual group $\widehat{G}$ is determined from the ($\widehat{G}$-invariant)
 Laurent monomials whose
exponents form the matrix ${}^t\!S$, and the dual polynomial 
$\widehat{W}$, which is a $\widehat{G}$-invariant polynomial in $m$ variables,
can be written
 explicitly as
\begin{equation}
\widehat{W}(y_0,\dots,y_{m-1})=\sum_{k=0}^{n-1}\widehat{c}_k
\prod_{j=0}^{m-1}y_j^{p_{jk}} \thinspace .
\end{equation}

This mirror partner is somewhat mysterious, due to the new parameters
$\widehat{c}_k$ which must be introduced.  However, the original group $G$
will act on those parameters (through its action on the set of mirror
polynomials), and the $G$-invariant quantities are familiar ones.  
Explicitly, if we write the moment map for the original $G$-action in the form
$\mu(x)=\frac12\sum_{k=1}^{n-1}\chi_k|x_k|^2$,
where $\chi_k$ is the character for the action of $G$ on the $k^{\text{th}}$
variable, then the invariant quantities for the $G$ action on the coefficients
of $\widehat{W}$ can be described as:
\begin{equation} \label{eq:invquant}
\frac1{2\pi i}\sum (\log \widehat{c}_k)\chi_k \in
\mathfrak{g}_{\mathbb{C}}^*/\mathfrak{g}_{\mathbb{Z}}^* \thinspace .
\end{equation}
(We have written the invariants additively, introducing a logarithm, and they
are thus multi-valued.)  Mirror symmetry predicts that
 the imaginary part of this invariant quantity \eqref{eq:invquant}
is to be identified with $r$, i.e.,
\begin{equation}
r=\frac{-1}{2\pi}\sum (\log |\widehat{c}_k|)\chi_k \thinspace .
\end{equation}
(Similarly, the invariant combinations $(\log z_\ell)/2\pi i$
of the original coefficients $c_j$ can
be identified with the complexification of the
K\"ahler parameters $\widehat{r}_\ell$ of the mirror
theory.)

This construction provides a global way to identify moduli spaces, and to
go beyond the small neighborhoods of large complex structure limit points.
The gauged linear sigma model makes sense for arbitrary values of $r$, not
just ones near an appropriate boundary point, and the description of a
mirror theory shows
that this realization could be a geometric one on the mirror partner.
In fact, an explicit (physics) computation can be made  \cite{MR96f:32036}
 of the locus where the theory becomes singular
(aside from toric  boundary points like $z=0$ and $z=\infty$),
and it reproduces the structure of the discriminant locus of the mirror
polynomial, including all of its components.
In the case of the quintic, there is only one component, a polynomial with
a single zero, at $z=-5^{-5}$.  It should be
stressed that this computation is made purely from the point of view
of the quintic itself, without reference to the mirror theory.

The quantum cohomology ring also corresponds as expected from mirror symmetry.  It can be
precisely calculated in the gauged linear sigma model on either side
(in one case from the data of the polynomial, in the other case from
the toric data, refined by analyzing the physics) and the results
agree \cite{MR96f:32036}.\footnote{A somewhat more physical argument for this
agreemenet of topological correlation functions 
was recently given in \cite{HoriVafa}, but we are  still lacking
is a direct argument that the mirror pair of gauged linear sigma
models give isomorphic physical theories, not just isomorphic
topological correlation functions \cite{MR98i:81236}.}
Relating this result to the enumerative predictions involves determining
an appropriate basis of cohomology (or in physical terms, calculating
the effect of renormalization), so one cannot derive the Mirror Theorem
directly in this way; however, the proofs of the Mirror Theorem rely
on similar results at some step along the way.

Applying this entire set-up to
 the case of the general quintic, we obtain  a $6\times126$
matrix;   the mirror partner can be determined from the transposed $126\times6$
matrix.  However, the calculation of the geometry of the mirror would be formidable
from this point of view.

An alternative is to begin with a quintic with fewer monomials.
If we start in $\mathbb{P}^4$ with the quintic defined by \eqref{eq:mirhyp},
then the associated factored matrix is given by 
\begin{equation}
\left(\begin{array}{cccccc}
1&1&1&1&1&1\\ 1&5&0&0&0&0\\ 1&0&5&0&0&0\\ 1&0&0&5&0&0\\
1&0&0&0&5&0\\ 1&0&0&0&0&5
\end{array}\right)=
\left(\begin{array}{ccccc}
1&1&1&1&1\\ 1&5&0&0&0\\ 1&0&5&0&0\\ 1&0&0&5&0\\ 1&0&0&0&5\\
1&0&0&0&0
\end{array}\right)
\left(\begin{array}{cccccc}
1&0&0&0&0&\hphantom{-}5\\ 0&1&0&0&0&-1\\ 0&0&1&0&0&-1\\
0&0&0&1&0&-1\\ 0&0&0&0&1&-1
\end{array}\right) \thinspace .
\end{equation}
(We can choose  $\mu=(1,0,0,0,0,0)$ and $\nu=(1,0,0,0,0,0)$ 
to obtain 
${}^t\!\nu P\mu=1$ and verify the conditions on the data.)
The group described by this factorization is $G=U(1)$.

To form the mirror, we take the transpose:
\begin{equation}
\left(\begin{array}{cccccc}
1&1&1&1&1&1\\ 1&5&0&0&0&0\\ 1&0&5&0&0&0\\ 1&0&0&5&0&0\\
1&0&0&0&5&0\\ 1&0&0&0&0&5
\end{array}\right)=
\left(\begin{array}{ccccc}
1&\hphantom{-}0&\hphantom{-}0&\hphantom{-}0&\hphantom{-}0\\ 
0&\hphantom{-}1&\hphantom{-}0&\hphantom{-}0&\hphantom{-}0\\ 
0&\hphantom{-}0&\hphantom{-}1&\hphantom{-}0&\hphantom{-}0\\ 
0&\hphantom{-}0&\hphantom{-}0&\hphantom{-}1&\hphantom{-}0\\ 
0&\hphantom{-}0&\hphantom{-}0&\hphantom{-}0&\hphantom{-}1\\
5&-1&-1&-1&-1
\end{array}\right)
\left(\begin{array}{cccccc}
1&1&1&1&1&1\\ 1&5&0&0&0&0\\ 1&0&5&0&0&0\\
1&0&0&5&0&0\\ 1&0&0&0&5&0
\end{array}\right) \thinspace .
\end{equation}
Once again, this represents the homogeneous polynomial  \eqref{eq:mirhyp}.
However, this time the group is $G=U(1)\times (\mathbb{Z}_5)^3$ and so
the Calabi--Yau is actually a hypersurface in a quotient space
$\mathbb{P}^4/(\mathbb{Z}_5)^3$.

\section{Consequences of the Homological Mirror Conjecture}

Kontsevich has proposed an intriguing extension of the original mirror
symmetry statements, comprising what he calls the ``homological mirror
conjecture'' \cite{MR97f:32040}. Briefly put, Kontsevich
employs the K\"ahler structure to define  Lagrangian
submanifolds of a Calabi--Yau manifold, and proposes that for a mirror
pair $(X,Y)$, Fukaya's $A_\infty$-category \cite{Fukaya} of Lagrangian
submanifolds of $X$ should be isomorphic to the bounded derived category of
coherent sheaves on $Y$.  This conjecture has proven to be very deep, with
results to date being primarily about the case of elliptic curves.
(See \cite{Manin} for a recent progress report.)

As the complex structure on $X$ is varied, the set of Lagrangian submanifolds
(with respect to a fixed K\"ahler structure) changes.  In fact, upon
traversing various loops in the complex structure moduli space, a large variety
of monodromy transformations on the set of Lagrangian submanifolds is
obtained.  In 1996, Kontsevich raised the question of reproducing this
effect on the other side of the proposed 
mirror symmetry relationship \cite{Kont-rutgers}: can we generate these
monodromy transformations as automorphisms of the bounded derived category of
coherent sheaves on $Y$?

One source of monodromy transformations we have already seen: each element
$\lambda\in H^2(X,\mathbb{Z})$ generates a monodromy transformation 
$\mathcal{T}_\lambda$
whose action on the even cohomology is given by
\begin{equation}
\mathcal{T}_\lambda:\gamma\mapsto \gamma\wedge e^{\lambda} \thinspace .
\end{equation}
(The exponential is actually a finite sum, since $\lambda^{\wedge k}=0$ 
for $k\gg0$.)

Kontsevich proposed another type of monodromy transformation.
Consider the diagonal $\Delta\subset X\times X$,
and its ideal sheaf $\mathcal{I}_{\Delta}$,
and consider the automorphism of the derived category defined by
\begin{equation}
\mathcal{E}\mapsto (p_2)_*\left(p_1^*\mathcal{E}\otimes 
\mathcal{I}_{\Delta}\right) \thinspace .
\end{equation}
The effect on cohomology, denoted by $\mathcal{S}$, can be described as
\begin{equation}
\mathcal{S}:\gamma\mapsto \gamma-\left(\int 
\gamma\wedge\operatorname{Todd}T_X\right)\cdot 1_X \thinspace ,
\end{equation}
where $1_X$ is the standard generator of $H^0(X,\mathbb{Q})$.

In the case of a quintic hypersurface, $H^{2k}(X,\mathbb{Q})$ is rank 
one with a standard
generator $\lambda^k$, for $k=0,1,2,3$.  With respect to this basis, 
we have matrices
\begin{equation}
{\renewcommand{\arraystretch}{1.1}
T_{\lambda} =
\left(\begin{array}{cccc}
1&1&\frac12&\frac16\\ &1&1&\frac12\\ & &1&1\\ & & &1\\
\end{array}\right),
\qquad
S =
\left(\begin{array}{cccc}
\hphantom{-}1& & & \\ -\frac{25}6&1& & \\ \hphantom{-}0& &1& \\ -5& & &1\\
\end{array}\right) \thinspace .
}
\end{equation}
Kontsevich then calculates the product of these matrices
\begin{equation}
{\renewcommand{\arraystretch}{1.1}
T_{\lambda} \cdot S = 
\left(\begin{array}{cccc}
-4&1&\frac12&\frac16\\ -\frac{20}3&1&1&\frac12\\ -5& 0&1&1\\ -5&0 &0 &1\\
\end{array}\right) \thinspace ,
}
\end{equation}
and observes that $(T_{\lambda}\cdot S)^5=I$.  

This is a really remarkable property!  It leads to the conclusion that in
matching these transformations to the monodromy on the complex structure
of the mirror (see 
Figure~\ref{fig:moduli}), the product of $T_\lambda$ and $S$ must
  represent the monodromy
around the Fermat point $z=\infty$.  Since  the matrix $T_\lambda$ 
is known to represent
the monodromy around the large complex structure limit point $z=0$,
$S$ must represent the monodromy around the point $z=-5^{-5}$. We have thus found
a natural way to describe this latter monodromy transformation.  This is 
a very
explicit indication of how mirror symmetry will reflect the structure
of the complex structure moduli space far from the small corners in which
mirror symmetry predictions are usually made.

More generally, Kontsevich suggests beginning
 with a holomorphic bundle $\mathcal{F}$ on $X$,
and defining an associated sheaf $\mathcal{I}_{\mathcal{F}}$ on $X\times X$ 
(with support on $\Delta$) as the kernel of the natural map
\begin{equation}
\mathcal{F}\boxtimes\mathcal{F}^*\to\mathcal{O}_{\Delta} \thinspace .
\end{equation}
There is then an endomorphism of the derived category defined by
\begin{equation}
\mathcal{E}\mapsto (p_2)_*\left(p_1^*\mathcal{E}\otimes 
\mathcal{I}_{\mathcal{F}}\right) \thinspace .
\end{equation}
This will actually be an automorphism provided that $\mathcal{F}$ satisfies
\begin{equation}
\operatorname{Ext}^k(\mathcal{F},\mathcal{F})=
\begin{cases}
\mathbb{C} & \text{if $k=0,n$ }\\
0 & \text{otherwise}
\end{cases} \thinspace .
\end{equation}
In particular, $\mathcal{F}=\mathcal{O}_X$ will give an automorphism of $X$
is $X$ is a ``proper'' Calabi--Yau manifold, i.e., one whose only convariantly
constant holomorphic forms are a $0$-form and an $n$-form.

The description sketched above has been made much more precise in
\cite{SeidelThomas,KhovanovSeidel}.  In particular, Seidel and
Thomas \cite{SeidelThomas} discovered a beautiful formula which shows
that the derived category admits actions of the braid group, which are
mirror to Dehn twists along Lagrangian spheres in the mirror manifold.
See \cite{Thomas} for a concise description of how this is related to
the mirror conjectures.

On the other hand, Horja \cite{Horja} generalized Kontsevich's picture
in a different way, and found a construction of automorphisms of the
derived category in the case of hypersurfaces in toric varieties which
is closely connected to the structure of the ambient toric variety.
Kontsevich's computation given above suggests that the case 
$\mathcal{F}=\mathcal{O}_X$ will generate an automorphism of the
derived category which is mirror to  the monodromy about a loop 
which surrounds the ``principal
component'' of the discriminant locus.  There are other components of
the discriminant locus, however, and each one is associated to a specific
subset $E$ of the toric variety $\mathbb{P}$, and to a mapping 
$\mathbb{P}\to\overline{\mathbb{P}}$
which
contracts $E$ to some subset $Z\subset\overline{\mathbb{P}}$.
(In the case of the ``principal component'', the subset is $E=\mathbb{P}$
and the contraction is the map from $\mathbb{P}$ to a point, with
$Z$=$\overline{\mathbb{P}}$=$\{\text{point}\}$.)

Horja shows that an appropriate sheaf on the fiber product
$E\times_Z E$ can be used to generate an automorphism of the derived
category of $X$.  In examples and special cases, he can check that this
coincides with the monodromy of periods on the mirror
around the corresponding component of the discriminant.

\section{Special Lagrangian Fibrations}

In 1996, a new geometric property of mirror pairs was deduced from physics by
Strominger, Yau, and Zaslow \cite{MR97j:32022}, through the study of so-called
{\it D-branes}.  This new property provides a purely geometric
way to define mirror pairs \cite{MR1714827}, at least in principle.  
Unfortunately,
our ability to calculate
with the geometric structures involved is very limited at present, and this
new ``geometric mirrror symmetry'' has not yet been fully connected up 
with other
constructions.  Geometric mirror symmetry has been analyzed in detail
for the ``Borcea--Voisin
threefolds'' \cite{MR99k:14066}, and some progress has been made
\cite{MR2000e:14066,GrossSLAGII,MR99f:81170,MR1738179,RuanI,%
GrossMirror,RuanIII,Ruanhyp}
in understanding the relationship between geometric mirror symmetry and
 Batyrev's mirror symmetry for toric hypersurfaces
or the gauged linear sigma model. 

The Strominger--Yau--Zaslow property is formulated in terms of the
``special Lagrangian submanifolds'' of a Calabi--Yau manifold $X$.  Every
Calabi--Yau manifold of dimension $2n$, when equipped with a Calabi--Yau
metric, admits both a covariantly constant K\"ahler form $\omega$ (the
existence of which reduces the holonomy from $SO(2n)$ to $U(n)$) and a
covariantly constant complex $n$-form $\Omega$ (the existence of which
further reduces the holonomy to $SU(n)$).  A {\it special Lagrangian
submanifold}\/ $M\subset X$ is a Lagrangian submanifold with respect to
$\omega$ (i.e., $\omega|_M\equiv0$) such that
$\operatorname{Im}(e^{i\theta}\Omega)$ vanishes on $M$---or equivalently,
that $\operatorname{Re}(e^{i\theta}\Omega)$ restricts to a constant
multiple of the volume form on $M$---for an appropriately chosen $\theta$.

Sadly, very few explicit examples of (compact) special Lagrangian
submanifolds of compact Calabi--Yau manifolds are known.  One construction
\cite{Br:min} represents the submanifold as the set of real points on an
algebraic variety defined over $\mathbb{R}$ (but these are rather rare),
while another only applies to the case of K3 surfaces \cite{HL} (or more
generally hyper-K\"ahler manifolds \cite{MR2000c:32075}).

In spite of our current lack of tools for constructing and analyzing
special Lagrangian submanifolds, we can still make the following
definition.  A {\it geometric mirror pair}\/ is a pair $(X,Y)$ of manifolds
equipped with Calabi--Yau metrics for which there are proper maps $f:X\to
B$ and $g:Y\to B$ to an $n$-manifold $B$, the general fibers of which are
special Lagrangian $n$-tori; moreover, the $n$-tori $f^{-1}(b)$ and
$g^{-1}(b)$ for generic $b$ should be {\it $T$-dual}\/ in the sense
explained in Section~\ref{sec:duality}.  (That is, their K\"ahler classes
should be related by the transformation which---on an appropriate sublocus
of the space of metrics---is described by $(r_1,\dots,r_n)\mapsto
(1/r_1,\dots,1/r_n)$.)  This definition makes no mention of $B$-fields, but
they can also be included in a natural way \cite{MR1714827}.

In the case of elliptic curves ($n=1$), the metric is flat and a special
Lagrangian $S^1$ is just a closed geodesic.  As is well known, if the
homology class is fixed then there is a fibration of the elliptic curve
over $S^1$ by closed geodesics in the specified class, with no singular
fibers.

In the case of K3 surfaces ($n=2$), every special Lagrangian
$T^2$-fibration can be interpreted \cite{HL} as a holomorphic elliptic
fibration for an approriate complex structure on the K3 surface, compatible
with the given Ricci-flat metric.  (In fact, any generic Ricci-flat metric
on a K3 surface admits such a fibration \cite{MR1714827}.)  In this case,
thanks to work of Kodaira \cite{MR32:1730}, a complete classification of
possible singular fibers is known: they are characterized by the conjugacy
class of the monodromy action on $H^1(T^2,\mathbb{Z})$.  The simplest
fibers, called {\it semistable}, are associated to unipotent monodromy
transformations.  In an appropriate basis, the monodromy matrix takes the
form
\begin{equation}
M=\left(\begin{array}{cc} 1&k\\0&1 \end{array} \right) \thinspace .
\end{equation}
For the generic elliptic fibration on a K3 surface, there are exactly $24$
semistable fibers, each with $k=1$.

The mirror partner of a given K3 surface (with a fixed Ricci-flat metric)
is another K3 surface with a different Ricci-flat metric.  The monodromy
matrices $M$ are replaced by ${}^t\!M^{-1}$; since this is congruent to
$M$, the monodromy data does not change.  In fact, it can be checked
\cite{MR1714827} (appealing to some work of Mukai \cite{Muk:bundles}) that
this mapping of geometric mirror partners is precisely the same one
mentioned in Section~\ref{sec:duality}, whose existence was originally
inferred from other considerations in physics.

In higher dimension, we are severely hampered at present by our lack of
tools for constructing and manipulating special Lagrangian submanifolds.
However, assuming the existence of a special Lagrangian $T^n$-fibration,
and assuming that for generic choices of Calabi--Yau metric the fibration
will have very generic monodromy behavior, it is possible to analyze that
behavior and obtain a clear picture of the topology of the fibration.  This
program has been carried out explicitly for the quintic threefold
\cite{MR1738179,RuanI,GrossMirror,RuanIII,Ruanhyp}.

Generic topological $T^3$-fibrations on Calabi--Yau threefolds take the
following form.\footnote{We are describing the analogue of $T^2$-fibrations
of a K3 surfaces which have exactly 24 singular fibers of the simplest
type.}  The base of the fibration $B$ is a $3$-sphere, and the discriminant
locus (where the $T^3$ fibers become singular) forms a graph $\Gamma$ on
$B$.  The edges of $\Gamma$ correspond to degenerations which do not affect
some particular $S^1$ within $T^3$, and look like the product of the $k=1$
degeneration of elliptic curves with a cylinder (the cylinder being the
product of the given edge and the unaffected $S^1$).  In particular, each
of the monodromy matrices is conjugate to
\begin{equation}
M=\left(\begin{array}{ccc} 1&1&0\\0&1&0\\0&0&1 \end{array} \right) \thinspace.
\end{equation}

The vertices of $\Gamma$ are all trivalent, and there are two types of
vertex.  The ``type $(2,1)$'' or ``type II'' vertices have the property
that the monodromy actions on $H^1$ associated to the three edges meeting
at the vertex have a common $2$-dimensional fixed plane, while the
monodromy actions on $H^2$ have fixed planes whose intersection is
$1$-dimensional.  For the ``type $(1,2)$'' or ``type III'' vertices this is
reversed: the monodromy actions on $H^1$ have fixed planes whose
intersection is $1$-dimensional, while the monodromy actions on $H^2$ have
a common $2$-dimensional fixed plane.

\begin{figure}
\begin{center}
\setlength{\unitlength}{.38em}

\begin{picture}(40,30)

\thinlines

\put(4,-2){\circle{1}}
\put(12,-2){\circle{1}}
\put(20,-2){\circle{1}}
\put(28,-2){\circle{1}}
\put(36,-2){\circle{1}}

\put(4,2){\circle*{1}}
\put(12,2){\circle*{1}}
\put(20,2){\circle*{1}}
\put(28,2){\circle*{1}}
\put(36,2){\circle*{1}}

\put(4,2){\line(0,-1){4}}
\put(12,2){\line(0,-1){4}}
\put(20,2){\line(0,-1){4}}
\put(28,2){\line(0,-1){4}}
\put(36,2){\line(0,-1){4}}

\put(0,4){\circle{1}}
\put(8,4){\circle*{1}}
\put(16,4){\circle*{1}}
\put(24,4){\circle*{1}}
\put(32,4){\circle*{1}}
\put(40,4){\circle{1}}

\put(40,4){\line(0,1){3}}
\put(40,4){\line(1,0){3}}
\put(36,10){\line(0,1){3}}
\put(36,10){\line(1,0){3}}
\put(32,16){\line(0,1){3}}
\put(32,16){\line(1,0){3}}
\put(28,22){\line(0,1){3}}
\put(28,22){\line(1,0){3}}
\put(24,28){\line(0,1){3}}
\put(24,28){\line(1,0){3}}

\put(0,4){\line(0,1){3}}
\put(0,4){\line(-1,0){3}}
\put(4,10){\line(0,1){3}}
\put(4,10){\line(-1,0){3}}
\put(8,16){\line(0,1){3}}
\put(8,16){\line(-1,0){3}}
\put(12,22){\line(0,1){3}}
\put(12,22){\line(-1,0){3}}
\put(16,28){\line(0,1){3}}
\put(16,28){\line(-1,0){3}}

\put(0,4){\line(2,-1){4}}
\put(4,2){\line(2,1){4}}
\put(8,4){\line(2,-1){4}}
\put(12,2){\line(2,1){4}}
\put(16,4){\line(2,-1){4}}
\put(20,2){\line(2,1){4}}
\put(24,4){\line(2,-1){4}}
\put(28,2){\line(2,1){4}}
\put(32,4){\line(2,-1){4}}
\put(36,2){\line(2,1){4}}

\put(8,8){\circle*{1}}
\put(16,8){\circle*{1}}
\put(24,8){\circle*{1}}
\put(32,8){\circle*{1}}

\put(8,8){\line(0,-1){4}}
\put(16,8){\line(0,-1){4}}
\put(24,8){\line(0,-1){4}}
\put(32,8){\line(0,-1){4}}

\put(4,10){\circle{1}}
\put(12,10){\circle*{1}}
\put(20,10){\circle*{1}}
\put(28,10){\circle*{1}}
\put(36,10){\circle{1}}

\put(4,10){\line(2,-1){4}}
\put(8,8){\line(2,1){4}}
\put(12,10){\line(2,-1){4}}
\put(16,8){\line(2,1){4}}
\put(20,10){\line(2,-1){4}}
\put(24,8){\line(2,1){4}}
\put(28,10){\line(2,-1){4}}
\put(32,8){\line(2,1){4}}

\put(12,14){\circle*{1}}
\put(20,14){\circle*{1}}
\put(28,14){\circle*{1}}

\put(12,14){\line(0,-1){4}}
\put(20,14){\line(0,-1){4}}
\put(28,14){\line(0,-1){4}}

\put(8,16){\circle{1}}
\put(16,16){\circle*{1}}
\put(24,16){\circle*{1}}
\put(32,16){\circle{1}}

\put(8,16){\line(2,-1){4}}
\put(12,14){\line(2,1){4}}
\put(16,16){\line(2,-1){4}}
\put(20,14){\line(2,1){4}}
\put(24,16){\line(2,-1){4}}
\put(28,14){\line(2,1){4}}

\put(16,20){\circle*{1}}
\put(24,20){\circle*{1}}

\put(16,20){\line(0,-1){4}}
\put(24,20){\line(0,-1){4}}

\put(12,22){\circle{1}}
\put(20,22){\circle*{1}}
\put(28,22){\circle{1}}

\put(12,22){\line(2,-1){4}}
\put(16,20){\line(2,1){4}}
\put(20,22){\line(2,-1){4}}
\put(24,20){\line(2,1){4}}

\put(20,26){\circle*{1}}

\put(20,26){\line(0,-1){4}}

\put(16,28){\circle{1}}
\put(24,28){\circle{1}}

\put(16,28){\line(2,-1){4}}
\put(20,26){\line(2,1){4}}

\end{picture}
\end{center}
\caption{Part of a graph for a generic $T^3$-fibration on a quintic threefold.}\label{fig2}
\end{figure}
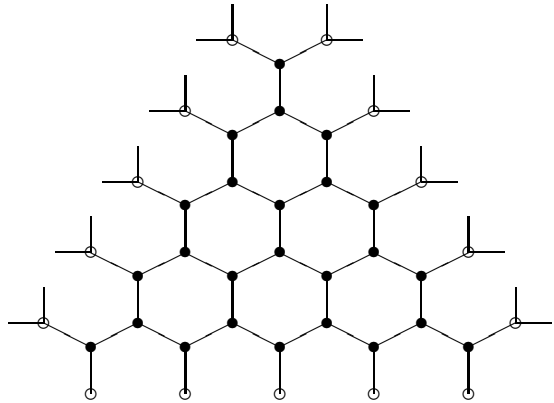

The graphs associated to generic $T^3$-fibrations on a quintic threefold
contain 250 vertices of type $(2,1)$ and 50 vertices of type $(1,2)$,
joined by a total of 450 edges.  A portion of one such graph, showing 25
vertices of the type $(2,1)$ vertices and 15 of the type $(1,2)$ vertices,
is illustrated in Figure~\ref{fig2}.  
(This corresponds to the toric diagram on the right side of Figure~\ref{fig1}.)
We warn the reader that although the
portion we illustrate is a planar graph, the overall graph is {\em not}\/
planar.

The combinatorics of the graph depend on both a triagulation of the
original Gorenstein cone, and of its dual.  In particular, if we change the
triangulation we get a different graph.  Gross has analyzed this change
in detail \cite{GrossMirror}, and shows that the same combinatorics which
appears in analyzing a toric flop also governs the change of graph for
the $T^3$-fibration.

Mirror symmetry replaces each monodromy transformation $M$ with the transformation
${}^t\!M^{-1}$.  This has the  geometric effect of exchanging the two
types of vertex, replacing a type $(2,1)$ vertex with one of type $(1,2)$,
and vice versa.  Mirror symmetry is thus realized in a simple and pleasing
geometric fashion.

It remains an interesting and important task for the future to relate this picture
of mirror symmetry to those which have come before.


\ifx\undefined\bysame
\newcommand{\bysame}{\leavevmode\hbox to3em{\hrulefill}\,}
\fi

\end{document}